\documentclass[12pt,a4paper,english]{smfart} % Font sizes: 10, 11, or 12; paper sizes: a4paper, letterpaper, a5paper, legalpaper, executivepaper or landscape; font families: sans or roman

\usepackage{fullpage}

\usepackage[table]{xcolor}

\usepackage[utf8]{inputenc}

\usepackage[OT2,T1]{fontenc}

\usepackage{amsfonts,lmodern}
\usepackage{amssymb,mathabx}
\usepackage{ifthen}
\usepackage{pstricks}
\usepackage{latexsym,eurosym,ulem}
\usepackage[french]{babel} 
\usepackage{footmisc}

\usepackage{alltt}
\usepackage{calrsfs}
\usepackage{fancyhdr}
\usepackage{graphicx}

\usepackage{amssymb}
\usepackage{amsmath}
\usepackage{smfthm,mathabx}
\usepackage{enumerate}
\usepackage{amsthm}
\usepackage{amsfonts}
\usepackage{graphicx}
\usepackage[colorlinks=true,linkcolor=black,citecolor=black,urlcolor=black]{hyperref}
\usepackage{fancyhdr}
\usepackage[all,cmtip]{xy}
\usepackage{alltt}
\usepackage{footmisc}
\usepackage{smfthm}

\font\rurm=wncyr10 scaled \magstep1
\def\sha{{{\textnormal{\rurm{Sh}}}}}

\def\Sha{{\sha}^2}
\def\Q{{\set Q}}

\def\F{{\set F}}

\def\k{{\rm k}}
\def\K{{\rm K}}
\def\L{{\rm L}}
\def\M{{\rm M}}
\def\KK{{\rm F}}
\def\HH{{\Gamma}}

\def\Z{{   \mathbb Z }}

\def\Zp{ {\mathbb Z}/p}
\def\Q{{\mathbb Q}}

\def\F{{\mathbb F}}
\def\G{{\rm G}}

\def\emp{{\emptyset}}

\def\J{{\rm J}}

\def\O{{\mathcal O}}

\def\Gal{{\mathrm{Gal}}}
\def\Ker{{\mathrm{Ker}}}
\def\Ann{{\mathrm{Ann}}}
\def\Cl{{\mathrm{Cl}}}
\def\Gov{{\mathrm{Gov}}}
\def\Fr{{\mathrm{Fr}}}

\def\sha{{{\textnormal{\rurm{Sh}}}}}
\def\CyB{{{\textnormal{\rurm{B}}}}}

\def\Sha{{\sha}^2}

\newtheorem{Theorem}{Theorem}

\newtheorem{Definition}{Definition}
\newtheorem{Fact}{Fact}
\newtheorem*{IntroTheo}{Theorem}
\newtheorem*{IntroCorollary}{Corollary}
\newtheorem*{Step}{Inductive Step}

\begin{document}

\date{\today}

\title{%On Ozaki's theorem on $p$-Hilbert class field towers\\ alternate: 
On Ozaki's theorem realizing prescribed $p$-groups %every %finite $p$-group 
as  $p$-class tower groups}
\thanks{
The second author  was partially supported by the ANR project FLAIR (ANR-17-CE40-0012) and  by the EIPHI Graduate School (ANR-17-EURE-0002). The third author was partially supported by Simons Collaboration grant \#524863. He also thanks FEMTO-ST for its hospitality and wonderful research environment during his visit there in the spring of 2022. All three authors were supported  for a Collaborate@ICERM visit in January, 2022. }
\author{Farshid Hajir, Christian Maire, Ravi Ramakrishna}
\address{Department of Mathematics \& Statistics, University of Massachusetts, Amherst, MA 01003, USA}
 \address{FEMTO-ST Institute, Universit\'e Bourgogne Franche-Comt\'e, CNRS,  15B avenue des Montboucons, 25000 Besan\c con, FRANCE} 
\address{Department of Mathematics, Cornell University, Ithaca, NY 14853-4201, USA}
\email{hajir@math.umass.edu, christian.maire@univ-fcomte.fr, ravi@math.cornell.edu}

\maketitle
\begin{abstract}
\black
We give a streamlined and effective proof of Ozaki's theorem that any finite $p$-group  $\HH$  is the Galois group of the $p$-Hilbert class field tower of some number field~$\KK$.  
Our work is inspired by Ozaki's and applies in broader circumstances. 
While his theorem is in the totally complex setting, we obtain the result in
any mixed signature setting for which there exists a number field $\k_0$ with class number prime to $p$. We construct
$\KK /\k_0$ by a sequence of $\Zp$-extensions ramified only at finite tame primes and
 also give explicit bounds on $[\KK:\k_0]$ and the number of ramified primes of $\KK/\k_0$ in terms of $\# \Gamma$.
 \end{abstract}

\section{Introduction}

For a number field $\k$, define $\L_p(\k)$ to be
 the compositum of all finite unramified  Galois $p$-extensions of $\k$. The extension $\L_p(\k)/\k$ is called the $p$-Hilbert class field tower of $\k$, and its Galois group $\Gal(\L_p(\k)/\k)$ is its $p$-class tower group. \black
%\blue
% Its maximal abelian quotient is isomorphic to the $p$-class group of $\k$. In analogy with the Cohen-Lenstra heuristics, for fields of fixed degree (or signature), one can study the distribution of $p$-class tower groups; for a survey of such results, see Wood \cite{Wood}. In such situations, there are severe restrictions on which groups can occur as $p$-class tower groups. 
In \cite{Ozaki}, Ozaki proved that 
%every finite p-group occurs as the p-class tower group of some totally complex number field F.
%showed that if the restriction on the degree of the base field is lifted, 
 every finite $p$-group   $\HH$ \black occurs as  the $p$-class tower group \black of some totally complex number field $\KK$. 
His strategy is as follows.
\black 
%\green DELETE ?   In Ozaki \cite{Ozaki},  every finite $p$-group   $\HH$  is realized as  the $p$-class tower group of some totally complex number field $\KK$. \black 
As finite $p$-groups are solvable, it is natural to proceed by induction.  After establishing the base case (realizing $\Z/p$ as a $p$-class tower group), it remains to show that given any short exact sequence of finite $p$-groups
\begin{equation} \label{eq:embedding}
1\to \Zp \to \G' \to \G \to 1
\end{equation}
where $\G:=\Gal(\L_p(\k)/\k)$, 
%and $\L_p(\k)$ is
%the compositum of all finite unramified  $p$-extensions of $\k$ for some number field $\k'$.
one can realize  $\G'$  as  $\Gal(\L_p(\k')/\k')$ for some number field $\k'$. Ozaki constructs such a 
 $\k'/\k$ via a sequence of carefully chosen $\Zp$-extensions.
 % ramified  at only finite tame primes. \green No ? he also uses wild ramification...
\black

\smallskip
 In this paper, we provide a streamlined and effective proof of Ozaki's theorem.  Some differences between our work and Ozaki's are: 
\begin{itemize}
\item   He must start with a totally complex $\k_0$ and then construct  a  field \black $\KK/\k_0$ whose $p$-Hilbert class field tower has the given $\HH$ as its Galois group, while we start with a number field $\k_0$ of arbitrary signature whose class number is prime to $p$. 
%\item  In \cite{Ozaki}, for $p$ odd, it is assumed $\mu_p \not \subset \k_0$. We allow $\mu_p \subset \k_0$ as long as  $(\#\C
\item Our result is effective and we are able %list all steps explicitly. In particular, this allows us 
to obtain  explicit upper bounds   on $[\KK:\k_0]$  and the number of ramified primes in  $\KK/\k_0$, all of which are tame and finite. 
\item  Moreover, we  bypass some  of the most  delicate and involved arguments of \cite{Ozaki}.%, we \red simplify the proof and reduce the size of $[\KK:\k_0]$. \black
\end{itemize}

We prove:

\begin{IntroTheo}
Let $\HH$ be a finite $p$-group and $\k_0$  a number field with $(\#\Cl_{\k_0},p)=1$. There exist infinitely many number fields ${\KK}/\k_0$ such that $\Gal(\L_p(\KK)/\KK) \simeq \HH$ and
\begin{itemize}
\item if $\mu_p\not \subset \k_0$ then $\KK/\k_0$ is of degree at most \black $p^2 \cdot \#\HH$ \black
and is ramified at at most \black $2+2\log_p(\#\HH)  $ \black finite tame primes,
\item if $\mu_p \subset \k_0$ then $\KK/\k_0$ is of degree at most \black $p \cdot (\#\HH)^2 $ \black
and is ramified at at most \black $1+3 \log_p(\#\HH)   $ \black finite tame primes.
\end{itemize}
\end{IntroTheo}

As any  countably generated  pro-$p$ group $\HH$ is the inverse limit of finite $p$-groups, Ozaki shows any such $\HH$ is the Galois group of the maximal unramified $p$-extension of some  infinite algebraic extension of $\Q$. The corresponding corollary of our theorem is: \black

\begin{IntroCorollary}\label{Corollary:totallyreal}
Any countably generated pro-$p$ group $\HH$, including $p$-adic analytic $\HH$, can be realized as the  $p$-class tower group \black of a totally real tamely ramified infinite extension $\KK/\Q$.
\end{IntroCorollary}
 
We now give  details about the structure of our proof and the difference between our methods and Ozaki's, though
we were very much inspired by Ozaki's beautiful theorem and techniques.

\smallskip
We start the base case of the inductive process with any number field $\k_0$, of any signature, whose class number is prime to $p$.  Referring to the group extension~\eqref{eq:embedding} with $\G$ being trivial,  one has to find an extension $\k'/\k_0$  such that $\k'$ has \black  $p$-class group tower exactly $\Zp$, which is equivalent to
the $p$-class  group being $\Zp$.
This is a standard argument and is part of Proposition~\ref{prop:finish}. 

\smallskip

The base case being done, we proceed to the inductive step (with our base field relabeled $\k$).  There are two cases, depending on whether~\eqref{eq:embedding} splits or not.   \black For the sake of brevity, we only \black outline the nonsplit case   in this introduction; the split case is handled similarly. For a set of places of $\k$, we say that an extension $\k'/\k$ is  exactly ramified at $S$ if it is ramified at all the places in $S$ and nowhere else.
\black
We need to find a suitable  tame prime $v_1$ of $\k$ such that
\begin{itemize}
\item  $v_1$ splits completely in $\L_p(\k)/\k$,
\item   There is no $\Zp$-extension of $\k$  exactly \black ramified at $v_1$,
\item  The maximal $p$-extension $\L_p(\k)_{ \{v_1\} }\black/\L_p(\k)$ exactly ramified at the primes of $\L_p(\k)$ above $v_1$ is of \black degree $p$ and solves the embedding problem~\eqref{eq:embedding}.
\end{itemize}

Arranging this and its split analog are the main technical difficulties. One then chooses a second prime $v_2$ that also solves the embedding problem as above and remains prime in $ \L_p(\k)_{ \{v_1\} }\black/\L_p(\k)$.  The existence of $v_1$ and $v_2$  will follow from  Chebotarev's Theorem.  The compositum of these two 
solutions, after a $\Zp$-base change  $\k'/\k$  ramified at both $v_1$ and $v_2$ (which exists!), gives the unramified solution to the embedding problem~\eqref{eq:embedding} which we show is $\L_p(\k')$.
This is done in the proof of Theorem~\ref{Theorem:Induction}.

\smallskip
 Our ability to choose primes $v_i$ as above depends upon the existence of {\it Minkowski units} in the tower $\L_p(\k)/\k$, namely  on the condition \color{black} that $\O^\times_{\L_p(\k)} \otimes \F_p \simeq \F_p[\G]^\lambda \oplus N$ where $N$ is an 
$\F_p[\G]$-torsion module and $\lambda$ is  a large enough integer. \black
In some situations, Minkowski units are rare - see \S 5.3 of \cite{HMRmink}. By contrast, both for Ozaki's proof (implicitly) and ours (explicitly), much of the work involves seeking fields for which they exist in abundance. \black
\black

\smallskip
 If $\mu_p \subset \k$, we may not be able to make our choices of $v_i$ as above to both split completely in $\L_p(\k)/\k$ and solve the  nonsplit  embedding problem~\eqref{eq:embedding}. In this case we need to perform an extra base change $\tilde{\k}/\k$ to shift  \black  the obstruction to the embedding problem so that we can proceed as above.
 The base change $\tilde{\k}/\k$ must preserve the tower, that is 
$\L_p(\tilde{\k})=\L_p(\k)\tilde{\k}$.  Theorem~\ref{Theorem:basechange} provides such a $\tilde{\k}$. \black

\smallskip
 Finally we check that the condition `$\lambda$ is large enough' persists, that is  there are enough Minkowski units to keep the induction going. Proposition~\ref{prop:finish} guarantees this.  To sum up, the key ingredients of the proof of the above Theorem and Corollary are Theorems~\ref{Theorem:basechange} and~\ref{Theorem:Induction} and Proposition~\ref{prop:finish}.
\black

\smallskip
We now explain in some detail Ozaki's approach and our simplifications.
\begin{itemize}
\item 
Using a result of Horie, \cite{Horie}, Ozaki starts with a quadratic imaginary field  with class number prime to $p$ in which $p$ is inert. He then chooses a suitable layer $\k$ in the cyclotomic $\Z_p$-extension  as the starting point of his induction. Assuming the problem solved for $\G$ in~\eqref{eq:embedding} and relabelling $\k$ as his base field, he proceeds inductively with the goal  to find a $\k'\supset \k$ whose
$p$-Hilbert class field tower has Galois group $\G'$. % as in~\eqref{eq:embedding}. 
For the induction to go forward, Ozaki needs $r_2(\k) \geq B_p(\k)$ 
 (implicit in this inequality is the existence of enough Minkowski units) \black
where $B_p(\k)$  is a certain explicit  quantity depending \black
 on $\k$, $\G$ and the $p$-part of the class group of $\K:=\L_p(\k)(\mu_p)$. \black
This involves delicate estimates in \S 4 of \cite{Ozaki}. 
We replace $r_2(\k)\geq B_p(\k)$  with $f(\k) \geq 2h^1(\G)+3$ where $h^i(\G):=\dim H^i(\G,\Zp)$  and $f(\k)$, which is a lower bound for the number of Minkowski units in  $\L_p(\k)/\k$, depends only on $h^1(\G)$, $h^2(\G)$ and the signature of $\k$.
We neither consider $\K$  nor invoke the estimates of \S 4 of \cite{Ozaki}.
\item In \S 5 of \cite{Ozaki}, Ozaki proves his base change Proposition $1$, namely he shows there exists a ramified 
$\Zp$-extension $\tilde{\k}/\k$ such that $\Gal(\L_p(\tilde{\k})/\tilde{\k}) \simeq \Gal(\L_p(\k)/\k)$. \black He uses this repeatedly when solving each embedding problem~\eqref{eq:embedding}. \black
Several tame primes are ramified in $\tilde{\k}/\k$ and he also needs that $\K$ and $\K\tilde{\k}$ have the same $p$-class group. This makes the proof significantly more involved. \black Theorem~\ref{Theorem:basechange} of this paper, our version of his Proposition $1$, has only one tame prime of ramification and  $\K$ plays no role. We only invoke  Theorem~\ref{Theorem:basechange} when $\mu_p \subset \k$. In particular, for $p$ odd, our Corollary above makes no use of Theorem~\ref{Theorem:basechange}.
\item To solve
the embedding problem~\eqref{eq:embedding}, 
Ozaki base changes several times (to a field relabeled $\k$) and then uses
 a  wildly ramified $\Zp$-extension $\L/\L_p(\k)$ 
 %with $\L/\k$ Galois, that is formed by essentially adjoining the  $p$th root of a unit of $\L_p(\k)$
  to solve~\eqref{eq:embedding}. After more base changes this is switched to a solution ramified at one tame prime. 
He then proceeds as in the description of this work using
two such solutions and a base change that absorbs the ramification at both tame primes
to find a $\k'$ such that $\Gal(\L_p(\k')/\k' )=\G'$. 
{\it We go directly to this last step and require at most two $\Zp$-base changes to solve the embedding problem. This allows us to quantify explicitly
both the degree and number of ramified primes of ${\KK}/\k$.}
\end{itemize}

\medskip
\black
 %\sout{It is worth remarking  that the existence of Minkowski units is not well-understood. For $p=2$ and quadratic imaginary  base \black fields they seldom exist - see \S 5.3 of \cite{HMRmink}. By contrast, both for Ozaki's proof and ours, they play a starring role and much of the work is to show that choices can be made to create situations where they occur with abundance. }
 %How one orders number number fields seems to affects their prevalence. 

\black

\medskip

{\bf Notations}

Let $p$ be a prime number.

$\bullet$ $\L$ is a number field, $\O_\L$ its ring of integers, $\O_\L^\times$ its units  and
$\Cl_\L$ is the $p$-Sylow of the class group of $\L$.

$\bullet$ For a finite set $S$ of primes of $\L$, set 
$$V_{\L,S}=\{x\in \L^\times, (x)={\mathcal I}^p, x\in (\L_v^\times)^p\,\, \forall v \in S\}.$$
In particular, one has the exact sequence:
$$1\longrightarrow \O_\L^\times\otimes \F_p \longrightarrow V_{\L,\emptyset}/(\L^\times)^p \longrightarrow \Cl_\L[p] \longrightarrow 1.$$

 $\bullet$ The superscript $^\wedge$ indicates the Kummer dual of an object $Z$ defined over a number field $\L$, though we never work with the $\Gal(\L(\mu_p)/\L)$ action on $Z^\wedge$.

\black 

 $\bullet$  $\L_S$  is the maximal pro-$p$-extension of $\L$ unramified outside $S$, $G_S :=\Gal(\L_S/\L)$ and 
 $\L_p(\L):= \L_\emp$, the maximal unramified pro-$p$-extension of $\L$, as it will ease notation at various points.

$\bullet$ $h^i(H):=\dim H^i(H,\Zp)$.

$\bullet$  $\Gov(\L):=\L(\mu_p)(\sqrt[p]{V_{\L,\emptyset}})$:  the governing field of $\L$. The span of 
$\{\Fr_v\}_{v\in S}$ in $M(\L):=\Gal(\Gov(\L)/\L(\mu_p))$ controls $\dim H^1(G_S)$.

\medskip 
The following may be helpful in orienting the reader:
\begin{itemize}
\item We frequently use finite tame primes with desired splitting properties in  number field extensions. We {\it always} use Chebotarev's theorem for the existence of such primes.
\item Our $\Zp$-extensions $\L'/\L$ of number fields are only ramified at (one or two) finite tame primes so $r_i(\L') =p \cdot r_i(\L)$ and 
$\mu_p \subset \L' \Longleftrightarrow \mu_p \subset \L$.
\item Note that $\k_0$ is our given base field, whereas $\k$ is a field used in the inductive process with $p$-class tower group $\G$ from~\eqref{eq:embedding}. Our task is to construct $\k'$ with $p$-class tower group $\G'$. Finally, $\tilde{\k}/\k$ is an extension having $p$-class tower group $\G$, the same as for $\k$.
\end{itemize}
\black
\section{Tools for the proof}

\black

\subsection{$\F_p[\G]$-modules and Minkowski Units}\label{subsection:MinkUs}

Let $\G$ be a finite group, a $p$-group in our situation. We record a few basic facts about finitely generated $\F_p[\G]$-modules~$M$. See \cite{CR}, \S $62$.
\begin{Fact} Any finitely generated $\F_p[\G]$-module $M$ is isomorphic to $\F_p[\G]^\lambda \oplus N$ where $N$ is a torsion $\F_p[\G]$-module   (every $n \in N$ is a torsion element)  and \black where $\lambda$ depends only on $M$.
\end{Fact}

Set $T_\G:=\sum_{g \in \G} g$.
Denote by $I_\G$ the augmentation ideal of $\F_p[\G]$. 
For $x\in M$ set $\Ann_\G(x):=\{ \alpha \in \F_p[\G] \mid \alpha \cdot x=0\}$.
Let $\{s_1,\cdots, s_{h^1(\G)}\}$ be a system of minimal generators of $\G$.  By Nakayama's lemma and  the fact that $I_\G/I_\G^2 \simeq \G/\G^p[\G,\G]$,  $I_\G$ can be generated, as $\G$-(right or left)-module,  by the elements $x_i:=s_i-1$.

\begin{prop} \label{prop_norm} With the $x_i$ as above,
let $M=\F_p[\G]^{h^1(\G)}$ %module on   generators  $e_1,\cdots, e_{h^1(\G)}$.
and $x= (x_1, x_2,\cdots, x_{h^1(\G)}) \in M$. 
 Then $\Ann_\G(x)=\F_p T_\G$. 
\end{prop}
\begin{proof}
$\Ann_\G(x) = \displaystyle{\bigcap_{i} \Ann_\G(x_i)=\Ann_\G(\langle x_i\rangle^{h^1(\G)}_{i=1} )= \Ann_\G(I_\G)=\F_p T_\G}$.
\end{proof}

 \begin{prop}\label{prop:prop_ann}  Let $M= \F_p[\G]^\lambda \oplus N$ be a finitely generated $\F_p[\G]$-module where $N$ is  torsion. Then $T_\G(M)  \simeq  \F_p^\lambda$.

 \end{prop}

 \begin{proof}  It is clear that $T_\G(\F_p[\G]^\lambda) \simeq  \F_p^\lambda$.  
 \black We now show $T_\G(N)=0$. 

Let $n \in N$ so  $\Ann_\G(n) \neq 0$. Note that  $\Ann_\G(n) \subset \F_p[\G]$ is a $p$-group stable under the action of 
the $p$-group $\G$ and thus has a fixed point. 
But it is easy to see the only fixed points of $\F_p[\G]$ are multiples of $T_\G$ so $T_\G \in \Ann_\G(n)$ as desired. 
 \end{proof}

\begin{Definition} We say  the tower $\L_p(\k)/\k$  with Galois group $\G$ has $\lambda$ Minkowski units if, as $\F_p[\G]$-modules, 
$V_{\L_p(\k),\emptyset}/\L_p(\k)^{\times p}  =  \O^\times_{\L_p(\k)}\otimes \F_p \simeq \F_p[\G]^\lambda \oplus N$ 
 where $N$ is an $\F_p[\G]$-torsion module.
\end{Definition}

\subsection{Extensions ramified at a tame set of primes}\label{subsection:tame}
We recall a standard formula  on the  number of $\Zp$-extensions of a number field with given tame ramification.
See \S $11.3$ of \cite{Koch} for a proof. Recall that for a field $\L$,
$\delta(\L)=\left\{\begin{array}{cc} 0 & \mu_p \not \subset \L\\1 &  \mu_p \subset \L \end{array}\right.$.

\begin{prop} \label{prop:h1dim}Let $\L$ be a number field, $p$ a prime number and $S$ a set of  tame primes of $\L$ prime to $p$. Then 
$$h^1(\G_{\L,S}):= \dim H^1(\G_{\L,S},\Zp) = \dim (V_{\L,S}/\L^{\times p})-r_1(\L)-r_2(\L) -\delta(\L) +1+\sum_{v\in S} \delta(\L_v).$$
\end{prop}

\black
Our $v\in S$  are always finite and have norm congruent to $1$ mod $p$ so $\delta(\L_v)=1$. 
Fact~\ref{Fact:GM} below follows immediately from Proposition~\ref{prop:h1dim} 
and the fact that  $\Gov(\L):=\L(\mu_p)(\sqrt[p]{ V_{\L,\emptyset}})$  is obtained by adjoining $p$th roots of elements of $\L$ to $\L(\mu_p)$.

\begin{Fact} \label{Fact:GM}     Let $S$ be a set of tame primes of $\L$ as above. For each $v\in S$ let $\Fr_v \in M(\L):=\Gal(\Gov(\L)/\L(\mu_p))$.  If the set
$\{\Fr_v, v\in S\}$ spans an $(\#S-d)$-dimensional subspace of $M(\L)$, then 
$$\dim H^1(\G_{\L,S},\Zp) =d+ \dim H^1(\G_{\L,\emp},\Zp).$$
When $\mu_p \not \subset \k$, $\Fr_v$ is only well-defined up to nonzero scalar multiplication.
\end{Fact}

\begin{Fact} \label{Fact:Chevalley} 
Let $\L$ be a number field such that $(\#\Cl_\L,p)=1$. 
Let $\L'/\L$ be a $\Z/p$-extension exactly ramified at $S=\{v_1,\cdots, v_r\}$ where the $v_i$ are finite and tame.
Then  $(\#\Cl_{\L'},p) =1$ if and only if  $\L'/\L$ is the  {\it unique } $\Z/p$-extension of $\L$ unramified outside $S$. In particular, that is the case when $|S|=1$.
\end{Fact}
\begin{proof}
\black Indeed, $(\#\Cl_{\L'},p) \neq 1$ if and only if there exists an unramified $\Z/p$-extension $H/\L'$ such that $H/\L$ is Galois (use the fact the the action of a $p$-group on a $p$-group always has  fixed points). Observe  that $H/\L$ cannot be cyclic of degree $p^2$ 
as all inertial elements of $\Gal(H/\L)$ have order $p$ and they would thus fix an unramified extension of $\L$, a contradiction.
So $\Gal(H/\L) \simeq \Z/p \times \Z/p$, and  $\L$ has at least two disjoints $\Z/p$-extension unramified outside $S$, also a contradiction. 
\end{proof}

Set $\CyB_{\L,S} = (V_{\L,S}/\L^{\times p})^{\wedge}$.   Recall 
$\sha^2_{\L,S} :=\Ker\left(H^2(G_S,\Zp) \to \oplus_{v\in S} H^2(G_v,\Zp)\right)$.
Fact~\ref{Fact:shafirst} below is well-known. See Theorem $11.3$ of \cite{Koch}. 
\black
\begin{Fact} \label{Fact:shafirst}
$\sha^2_{\L,S} \hookrightarrow \CyB_{\L,S}$.
\end{Fact}
Let $\lambda_\L$ be the number of Minkowski units in $\L_p(\L)/\L$.

\begin{Fact} \label{Fact:shasecond} 
If $\mu_p \not \subset \L$ then $\lambda_\L = r_1(\L)+r_2(\L)-1 +h^1(\G)-h^2(\G)$.\\
If $\mu_p \subset \L$ then $\lambda_\L \geq r_1(\L)+r_2(\L)-h^2(\G)$.
\end{Fact}

This result is Theorem $2.9$ of  \cite{HMRmink}, but we sketch the proof for the sake of keeping this paper self-contained.

\begin{proof}
 Set $\G=\Gal(\L_p(\L)/\L)$. We consider two ``norm maps'' induced by the norm map on units:
 $\O_{\L_p(\L)}^\times \to \O_{\L}^\times$.
 \begin{enumerate}
  \item[$-$] $ N_\G $ sending $\O_{\L_p(\L)}^\times\otimes \F_p $ to 
 $\displaystyle \frac{\O_\L^\times}{\O_\L^\times \cap (\O_{\L_p(\L)}^{\times})^p }\subset \O_{\L_p(\L)}^\times\otimes \F_p $;  
 \item[$-$] $ N_\G': \O_{\L_p(\L)}^\times\otimes \F_p  \rightarrow \O_{\L}^\times \otimes \F_p $. 
 \end{enumerate}
One easily sees   $ N_\G'  (\O_{\L_p(\L)}^\times \otimes \F_p) \twoheadrightarrow 
 N_\G (\O_{\L_p(\L)}^\times \otimes \F_p)$ and this is an isomorphism provided  $\O_\L^\times \cap (\O_{\L_p(\L)}^{\times})^p =(\O_\L^{\times})^p$: in particular this is the case  when $\mu_p \not \subset \L$, see Proposition $2.8$ of \cite{HMRmink}.

 Write 
$\O_{\L_p(\L)}^\times \otimes \F_p \simeq \F_p[\G]^{\lambda_\L} \oplus N$, 
where $N $ is an  $\F_p[\G]$-torsion module. 
By Proposition~\ref{prop:prop_ann} one has  $ N_\G  (\O^\times_{\L_p(\L)} \otimes \F_p) \simeq  \F_p^{\lambda_\L}$.  
%(see the proof of Proposition \ref{prop:prop_ann}). 
Hence, when $\mu_p \not \subset \L$
$$\dim \left(  \frac{\O_\L^\times\otimes \F_p}{ N_\G'  (\O_{\L_p(\L)}^\times\otimes \F_p)}\right)=  \dim (\O_\L^\times\otimes \F_p)  - \lambda_\L.$$
When $\mu_p \subset \L$, note that the `difference' between the images of $N_\G$ and  $N_\G'$  has $p$-rank at most  
$\dim \left( \frac{\O_\L^\times\cap \O_{\L_p(\L)}^{\times p}}{(\O_\L^{\times})^p}\right) \leq h^1(\G)$, so 
$$ \dim\left( \frac{\O_\L^\times\otimes \F_p}{ N'_\G  (\O_{\L_p(\L)}^\times   ) }\right) \geq  \dim (\O_\L^\times\otimes \F_p) - \lambda_\L - h^1(\G).$$

To conclude, we use the well-known  equality  (see \cite[Lemma 9]{Roquette}): $$h^2(\G)-h^1(\G)= \dim \left(\frac{\O_\L^\times\otimes \F_p}{N'_\G (\O_{\L_p(\L)}^\times\otimes \F_p)}\right).$$
 \end{proof}

\subsection{Solving the  ramified  embedding problem with one tame prime}\label{subsection:embedding}

We  start with our nonsplit exact sequence:
\begin{eqnarray}\label{eqn:seq}
 1 \longrightarrow \Zp \longrightarrow \G' \longrightarrow \G \longrightarrow 1.
\end{eqnarray}
given by the element $0\neq \varepsilon \in H^2(\G,\Zp)$.

We assume that $\G=\Gal(\L_p(\k)/\k)$.

\smallskip

 Set  $S=\{v\}$ where $v$ is a finite tame prime of $\k$. 
We first show the existence of a lift of $\G$ to $\G'$ in some $\k_S/\k$ for certain  $v$ of $\k$. We call this solving the embedding problem \eqref{eqn:seq} in $\k_S$.

Recall that $\Sha_{\k,S} \hookrightarrow \CyB_{\k,S}$ by Fact~\ref{Fact:shafirst}. Here 
$\Sha_{\k,\emptyset} \simeq H^2(\G_{\k,\emptyset},\Zp)\simeq H^2(\G,\Zp)$. 
Let $Inf_S: H^2(\G_{\k,\emptyset},\Zp) \rightarrow H^2(\G_{\k,S},\Zp)$ be the inflation map.
We  have the commutative diagram:

$$ \xymatrix{& \Sha_{\k,\emp} \ar@{->}[r]^{Inf_S} \ar@{^{(}->}[d]^{h} &\Sha_{\k,S} \ar@{^{(}->}[d]^{g}\\
 \big (\k_v^\times \otimes \F_p\big)^\wedge \ar@{->}[r]& \CyB_{\k ,\emptyset}  \ar@{->>}[r]^{f_S}&  \CyB_{\k,S}
}$$

By Hoeschmann's criteria   (see \cite[Chapter 3, \S 5]{NSW}), \black the embedding problem has a solution in $\k_S$ if and only if $Inf_S(\varepsilon) = 0$.
As $\L_p(\k)/\k$ is unramified, $Inf_S(\varepsilon) \in \Sha_{\k,S}$ and as
$g(Inf_S(\varepsilon))=f_S(h(\varepsilon))\in \CyB_{\k,S}$, the embedding problem has a solution if and only if $h(\varepsilon) \in \Ker(f_S)$.

 Set $\Gov_S(\k):=\k(\mu_p)(\sqrt[p]{ V_{\k,S}})$.
In the governing extensions $\k(\mu_p)\subset \Gov_S(\k) \subset \Gov(\k)$, one sees that \black the kernel of the map $f_S:\CyB_{\k,\emp} \twoheadrightarrow \CyB_{\k,S}$ is exactly the (unramified) decomposition group $D_v$  of the 
 prime $v$. As noted in Fact~\ref{Fact:GM}, if $w_1,w_2|v$ are two primes of $\k(\mu_p)$, their Frobenii in $\Gal(\Gov(\k)/\k(\mu_p))$ differ by a nonzero scalar multiple.

We have proved 
\begin{lemm} \label{lemm:lemm_lifting}
 The embedding problem \eqref{eqn:seq} has a solution in $\k_S/\k$ if and only if  $h(\varepsilon) \in D_v$. 
 Thus it has a solution in $\k_S/\k$ if we choose the prime $v$ such that $\langle \Fr_v  \rangle=\langle h(\varepsilon) \rangle$ in $M(\k)$,  that is the lines spanned by these elements in $M(k)$ are equal. This is always possible by Chebotarev's Theorem. 
\end{lemm}

\subsection{Cohomological facts implying the persistence of Minkowski units}
Our main aim in this paper is to show that given a short exact sequence 
 $$1 \to \Zp \to \G' \to \G \to 1$$ of finite $p$-groups where 
 $\G =\Gal(\L_p(\k)/\k)$,  there exists a finite tamely ramified extension $\k'/\k$ with
 $\G' = \Gal(\L_p(\k')/\k')$. 
 To solve this embedding problem using Theorem~\ref{Theorem:Induction},  the tower $\L_p(\k)/\k$ must have $2h^1(\G)$ Minkowski units.
 Proposition~\ref{prop:growth} below shows that if we start with enough Minkowski units, after a 
 base change that realizes $\G'$, we will be able to continue the induction.  Proposition~\ref{prop:stability}, which is only needed in the case when $\mu_p \subset \k$,  shows that given at least $h^1(\G)$ Minkowski units,  we can perform a base change that preserves the tower and the number of Minkowski units increases. 
 Proposition~\ref{prop:grouptheory}
 is a basic group theory result  bounding  $h^1(\G')$ and $h^2(\G')$ in terms  of $h^1(\G)$ and $h^2(\G)$.
 Furata proves a similar result in Lemma $2$ of \cite{Furata}. 

\smallskip
%Recall $h^i(H):= \dim H^i(H,\Zp)$ and 
Set
$H^2(\G',\Zp)_1:=\Ker\left(H^2(\G',\Zp) \stackrel{Res}{\to} H^2(\Zp,\Zp)\right)$. Note $h^2(\Zp)=1$  so
$h^2(\G')_1$ is either $h^2(\G')$ or $h^2(\G')-1$ and in either case $h^2(\G')_1 \geq h^2(\G')-1$.

\begin{prop} \label{prop:grouptheory}
Let $$1 \to \Zp \to \G' \to \G \to 1$$ be a short exact sequence of finite $p$-groups. Then
$h^1(\G') \leq h^1(\G)+1$ and
$h^2(\G') \leq h^1(\G)+h^2(\G)+1$.
\end{prop}
\begin{proof} 
The $h^1$ result is clear. For the $h^2$ statement
we have the long exact sequence 
(see for instance \cite{DHW})
$$0 \to H^1(\G,\Zp) \to H^1(\G',\Zp)\to H^1(\Zp,\Zp)^G $$
$$\to H^2(\G,\Zp)\to H^2(\G',\Zp)_1 \to H^1(\G,H^1(\Zp,\Zp)).$$
\black If $\G'\to \G$ splits,  we have \black 
$$ 0 \to H^2(\G,\Zp)\to H^2(\G',\Zp)_1 \to H^1(\G,H^1(\Zp,\Zp))$$
so $h^2(\G')_1\leq h^2(\G)+h^1(\G)$ and since $h^2(\G')_1 \geq h^2(\G')-1$ the result follows.

In the nonsplit case we have
$$0 \to H^1(\Zp,\Zp)^\G \to H^2(\G,\Zp)\to H^2(\G',\Zp)_1 \to H^1(\G,H^1(\Zp,\Zp))$$
so $h^2(\G')_1 \leq h^2(\G)-1+ h^1(\G)$ so
$h^2(\G') \leq h^1(\G)+h^2(\G)$.
\end{proof}

\begin{Definition}\label{Definition:f}
For a number field $\L$ set $\G=\Gal(\L_p(\L)/\L)$. Define  $f$   as follows:
$$f(\L)=\left\{\begin{array}{lc} r_1(\L)+r_2(\L)-h^2(\G)+h^1(\G) -1& \mu_p\not \subset \L\\
r_1(\L)+r_2(\L)-h^2(\G) & \mu_p\subset \L \end{array}\right.  .$$ Fact~\ref{Fact:shasecond} 
implies $f(\L)$ is a lower bound on the number of Minkowski units of    $\L_p(\L)/\L$.
\end{Definition}

\begin{prop}\label{prop:stability} Let $ \tilde{\k} \black /\k$ be a  $\Zp$-extension ramified at finite tame primes such that $\G=\Gal(\L_p(\k)/\k) = \Gal(\L_p(\tilde{\k})/\tilde{\k})$. 
Then $f(\tilde{\k}) = f(\k)+(p-1)(r_1(\k)+r_2(\k))$.
\end{prop}
\begin{proof} This follows immediately as we have the same group $\G$ for $\k$ and $\tilde{\k}$, $\mu_p\subset \tilde{\k} \iff \mu_p \subset \k$ and $r_i(\tilde{\k})=p\cdot r_i(\k)$.
\end{proof} 

\begin{prop}\label{prop:growth} Let $\k'/\k$ be a  tamely ramified $\Zp$-extension such that $\G=\Gal(\L_p(\k)/\k)$ and $\G'= \Gal(\L_p(\k')/\k')$
where 
$$1\to \Zp\to \G'\to \G\to 1.$$ 
Let $f(\k)$ be as in Definition~\ref{Definition:f}.
Then $$f(\k) \geq 2h^1(\G)+3 \implies f(\k')\geq 2h^1(\G')+3.$$
\end{prop}
\begin{proof} 
We do the case  $\mu_p \not \subset \k$ first.
We need to prove
$$r_1(\k)+r_2(\k)-h^2(\G)+h^1(\G) -1\geq 2h^1(\G)+3$$ 
$$ \implies  r_1(\k')+r_2(\k')-h^2(\G')+h^1(\G') -1\geq 2h^1(\G')+3,$$ 
that is
$$ r_1(\k')+r_2(\k')\stackrel{?}{\geq} h^1(\G')+h^2(\G')+4. $$
\black Clearly
$$r_1(\k')+r_2(\k')=  p(r_1(\k)+r_2(\k)) \geq p(h^1(\G)+h^2(\G)+4)$$ and
by Proposition~\ref{prop:grouptheory} we have
$$h^2(\G')+h^1(\G')+4 \leq (h^1(\G)+h^2(\G)+1)+(h^1(\G)+1)+4=2h^1(\G)+h^2(\G)+6$$
so it suffices to show 
$$(p-1)h^2(\G)+(p-2)h^1(\G) +4p\stackrel{\checkmark}{\geq} 6.$$
\black
This holds for all $p$.

\vskip1em
When $\mu_p \subset \k$.
We need to prove
$$r_1(\k)+r_2(\k)-h^2(\G)\geq 2h^1(\G)+3 \implies  r_1(\k')+r_2(\k')-h^2(\G')\geq 2h^1(\G')+3,$$ 
that is
$$ r_1(\k')+r_2(\k')\stackrel{?}{\geq} 2h^1(\G')+h^2(\G')+3. $$
Again using  Proposition~\ref{prop:grouptheory} and that $r_i(\k')=p\cdot r_i(\k)$ 
 it suffices to show 
$$(p-1)h^2(\G)+(2p-3)h^1(\G) +3p\stackrel{\checkmark}{\geq} 6$$
which holds for all $p$.
\end{proof} 

Proposition~\ref{prop:finish} below provides the base case of the induction.

\begin{prop}\label{prop:finish} Recall $(\# \Cl_{ \k_0},p)=1$.
There exists a tamely  ramified extension $\k'/\k_0$ such that 
\begin{itemize}
\item the $p$-part of the class group of $\k'$ is $\Zp$,
\item $[\k':\k_0]=p^3$,
\item and $f(\k') > 2h^1(\Zp)+3=5$.
\end{itemize}
\end{prop}
\begin{proof} Since $\L_p(\k_0)=\k_0$, we see $\G=\{e\}$. 
Choose a tame prime  $v$ of $\k$ whose Frobenius is trivial in the governing Galois group $M(\k)$.
By Fact~\ref{Fact:GM} there is  a unique $\Zp$-extension $\k_1/\k_0$. That
 $(\# \Cl_{\k_1},p)=1$  follows from Fact \ref{Fact:Chevalley}.  Repeat this process with $\k_1$ 
to get a field $\k_2$ with $(\#\Cl_{\k_2},p)=1$.

\smallskip
We  do one more base change to find a field $\k'$ with class group $\Zp$.
This is  proved more generally as part of Theorem~\ref{Theorem:Induction}, but we include a short proof here. 
 
 \smallskip
 Choose $v_1$ a finite tame prime of $\k_2$ with trivial Frobenius in $\Gov(\k_2)$ so that by Fact~\ref{Fact:GM} there exists a unique $D_1/\k_2$ ramified at $v_1$. As $D_1\cap \Gov(\k_2)=\k_2$, we may  choose $v_2$ a finite tame prime of $\k_2$ with trivial Frobenius in $\Gov(\k_2)$ such that 
 $v_2$ remains prime in $D_1/\k_2$. Again by Fact~\ref{Fact:GM} there exists a unique $D_2/\k_2$ ramified at $v_2$.

 \smallskip 
Let $D/\k_2$ be any of the $p-1$ `diagonal' $\Zp$-extensions of $\k_2$ between $D_1$ and 
$D_2$ %There are $p-1$ of these. 
so $D_1D_2/D$  is everywhere unramified. We claim $D_1D_2=\L_p(D)$. 
%This is a consequence of applying Fact \ref{Fact:Chevalley} twice.
  Indeed,
by Fact~\ref{Fact:Chevalley} applied to $D_1/\k_2$ we see $(\#\Cl_{D_1},p)=1$.
 As $v_2$ is inert in $D_1/\k_2$, the extension $D_2D_1/D_1$ is ramified only  at $v_2$ and  Fact~\ref{Fact:Chevalley} applied to $D_2D_1/D_1$  implies $(\#\Cl_{D_1D_2},p)=1$.
Whether or not $\mu_p\subset \k$, we have  ${\k'}:=D$, %with $[\k':\k]=p^3[\k:\mathbb Q]$  and 
$\Cl_{\k'} =\Zp$ 
and
$$f(\k') \geq r_1(\k') +r_2(\k') -h^2(\Zp) = p^3r_1(\k)+p^3r_2(\k) -1 \black >5=2h^1(\Zp)+3.$$ 
\end{proof}

Depending on $p$ and the signature of $\k_0$ one can decrease the number of base changes,   but this analysis complicates the statement of the main theorem without significant gain.

\section{Solving the embedding problem}\label{section:IS}

Having established the base case of our induction,   we now prove Theorem~\ref{Theorem:Induction}, the main

\begin{Step} Let $$1\to \Zp\to \G' \to \G \to 1$$ be exact and let $\k$ be a number field with $\Gal(\L_p(\k)/\k)=\G$ and
$f(\k) \geq 2h^1(\G)+3$. Then there exists a number field $\k'/\k$ with $\Gal(\L_p(\k')/\k')=\G'$ and
$f(\k') \geq 2h^1(\G')+3$. 
\end{Step}

Theorem~\ref{Theorem:basechange} below
is  only necessary for the key inductive step, Theorem~\ref{Theorem:Induction}, when $\mu_p \subset \k$.

\smallskip

Set $\K:=\L_p(\k)(\mu_p)$. We only consider finite tame primes $v$ of  $\k$ that split completely in $\K/\k$. 
%\smallskip
When $\mu_p \not \subset \k$, our  Frobenii in governing fields (or their subfields) are only defined  up to scalar multiples.
We write $\langle \Fr_{v} \rangle_{\Gov(\k)/\k(\mu_p)}$ for the well-defined line spanned by Frobenius at $v$ in
$\Gal(\Gov(\k)/\k(\mu_p))$. %It is well-defined. 
When the Frobenius is trivial there is no ambiguity so we write  $\langle \Fr_{v} \rangle_{\Gov(\k)/\k(\mu_p)}=0$.

\smallskip
We  need  primes $v$ of $\k$ that let us control $h^1(\Gal(\k_{ \{v\}}/\k))$ and
 $h^1(\Gal(\L_p(\k)_{ \{v\}}/\L_p(\k)))$ simultaneously
via Fact~\ref{Fact:GM}.
 Recall $M(\L_p(\k)):=\Gal(\Gov(\L_p(\k))/\L_p(\k)(\mu_p)) \simeq \F_p[\G]^{\lambda_\k} \oplus N$ where 
 $N$ is a torsion module over $\F_p[\G]$. We have no knowledge of $N$ and  must work with
the free part to control things over $\L_p(\k)$. 
%It is the presence of enough Minkowski units in $\L_p(\k)/\k$ that allows us  to satisfy the first point above over $\K$. 
We then use Proposition~\ref{prop:Ozaki} to control things over $\k$.

\black

\subsection{The Stability Theorem}
%Proposition~\ref{prop:Ozaki}  is only necessary for 
\begin{prop}\label{prop:Ozaki}
Let $F \subset \Gov(\L_p(\k))$ be the field fixed by $I_\G \cdot M(\L_p(\k))$.
%of the governing field of $\L_p(\k)$. 
%Then $\langle \Fr_v \rangle \subset M(\k)$ is well-defined
For $v$ of $\k$ splitting completely in $\K$ and $w|v$ in $\K$, the lines $\langle \Fr_w \rangle_{F/\K}$ do not dependent on  $w$ so we may write $\langle \Fr_v \rangle_{F/\K}$. Then 
$\langle \Fr_{v_1} \rangle_{F/\K}  =\langle \Fr_{v_2} \rangle_{F/\K}$ implies 
$\langle \Fr_{v_1} \rangle_{\Gov(\k)/\k(\mu_p)} = \langle \Fr_{v_2} \rangle_{\Gov(\k)/\k(\mu_p)} $.
If $\langle \Fr_{v_1} \rangle_{F/\K}  =0$ then $\langle \Fr_{v_1} \rangle_{\Gov(\k)/\k(\mu_p)} =0$.

\end{prop}

\begin{proof} 
This diagram is useful in Theorems~\ref{Theorem:basechange} and~\ref{Theorem:Induction} as well.
$${
\xymatrix{ 
& & & & \Gov(\L_p(\k)) &\\
& & & F\ar@{-}[ur] \ar@{-}@/^-1.50pc/[ur]_{I_\G \cdot \M(\L_p(\k))} \\
& &  \Gov(\k)\K\ar@{-}[d] \ar@{-}[ur]  &\\
& \K%:=\L_p(\k)(\mu_p)
\ar@{-}[ur] \ar@{-}@/^3.00pc/[uuurrr]^{\M(\L_p(\k))}& \Gov(\k) &\\
\L_p(\k)\ar@{-}[ur]^\Delta&  \k(\mu_p)\ar@{-}[ur] \ar@{-}[u]_\G&   &\\
\k  \ar@{-}[u]_\G \ar@{-}[ur]^\Delta&  &  &\\
%\k  \ar@{-}[uu]_\G \ar@{-}[ur]_\Delta&  & &\\
}}.$$ 

Let $\Delta =\Gal(\k(\mu_p)/\k) =\Gal(\K/\L_p(\k))$. As $\Gal(F/\K) := M(\L_p(\k))/I_\G \cdot M(\L_p(\k))$ is the maximal quotient of $M(\L_p(\k))$ on which $\G$ acts trivially, and 
%Fact~\ref{Fact:GM}  \green CHANGE/DELETE somethings ?? implies 
$\Delta$ acts on $\Gal(F/\K)$ by scalars, \black the line $\langle \Fr_{w} \rangle_{F/\K}$ is invariant under the action of $\Gal(\K/\k) =\G \times \Delta$. Since the $w|v$ form an orbit under this action of $\Gal(\K/\k)$,  this line is independent of the choice of $w|v$ as desired.

\smallskip
%Let $F$ be the fixed field of $I_\G\cdot  M(\L_p(\k)) \simeq I^\lambda_\G \oplus I_\G\cdot N \subset  M(\L_p(\k))$. 
As $\Gov(\k)\K/\K$ ascends from %the elementary $p$-abelian extension 
$\Gov(\k)/\k(\mu_p)$, 
%(it is possible,  when $\mu_p \subset \k$,  that $\K \cap \Gov(\k) \supsetneq \k(\mu_p)$ but that will play no role here),
we see $\G$ acts trivially on $ \Gal(\Gov(\k)\K/\K)$ 
%and as $\Gal(F/\K)$ is the maximal quotient of $M(\L_p(\k))$ on which $\G$ acts trivially, we see  
so $\Gov(\k)\K\subset F$. %By Proposition~\ref{prop:prop_ann}  
 Below, we  implicitly use that our primes of $\k$ split completely in $\K$. \black
If $\langle \Fr_{v_1} \rangle_{F/\K} =\langle \Fr_{v_2} \rangle_{F/\K} $,  these lines are equal when projected to
$\Gal(Gov(\k)\K/\K) \subset \Gal(\Gov(\k)\K/\k(\mu_p))$ and they are again equal in
$%M(\k):=
\Gal(\Gov(\k)/\k(\mu_p))$ so $\langle \Fr_{v_1} \rangle_{\Gov(\k)/\k(\mu_p)} =\langle \Fr_{v_2} \rangle_{\Gov(\k)/\k(\mu_p)} $. The last statement is clear. 
%The second part follows immediately as $v$ splits completely in $\Gov(\k)/\k$.
\end{proof}

\begin{Theorem}\label{Theorem:basechange} %With the hypotheses as in Proposition~\ref{prop:Ozaki} above,
Recall $\{x_i\}_{i=1}^{h^1(\G)}$ is a minimal set of generators of $I_\G$. 
Assume that 
$f(\k) \geq h^1(\G)$. 
Let $w$ be a degree one prime of $\K$ 
such that 
%$$\langle \Fr_{v} \rangle_{ Gov(\L_p(\k))/\K}  = \langle ((x_1,x_2,\cdots,x_{h^1(\G)},0,\cdots,0),0)\rangle_{ Gov(\L_p(\k))/\K}.$$
$$\Fr_w    = ((x_1,x_2,\cdots,x_{h^1(\G)},0,\cdots,0),0) \in M(\L_p(\k)) \simeq \F_p[\G]^{\lambda_\k}\oplus N.$$
Then for $v$ of $\k$ below $w$,
$\langle \Fr_v \rangle_{\Gov(\k)/\k(\mu_p)}=0$ so there exists a $\Zp$-extension $\tilde{\k}/\k$  ramified at $v$.   
%and $\L_p(\k)\tilde{\k}/\L_p(\k)$  is ramified   at  the $\F_p[\G]$ orbit of $w$,  that is all primes of $\L_p(\k)$ above $v$.
Furthermore, $\L_p(\tilde{\k})=\L_p(\k)\tilde{\k}$  and $f(\tilde{\k}) > f(\k)$. 
\end{Theorem}

\begin{proof} 
As $\Fr_{w} $ projects to $0$ in the $\F_p$-vector space $\Gal(F/\K)$,
Proposition~\ref{prop:Ozaki} implies
$\langle \Fr_v \rangle_{\Gov(\k)/\k(\mu_p)}=0$ so $\tilde{k}$ exists by Fact~\ref{Fact:GM}. We show  the $\F_p[\G]$-span of $(x_1,\cdots,x_{h^1(\G)}) \in\F_p[\G]^{h^1(\G)} $
%\subset \F_p[\G]^\lambda \oplus N$ 
has dimension $\#\G-1$ by computing the dimension of 
$ \cap^{h^1(\G)}_{i=1} \Ann(x_i)$. This intersection is the annihilator of $I_\G$ which 
 by Proposition \ref{prop_norm} is just $\F_pT_\G$, %    and   has $\F_p$-dimension one, 
 establishing our dimension result.
%, so the dimension of the $\F_p[\G]$-span of $(x_1,\cdots,x_{h^1(\G)})$ is $\# \G-1$. 
By Fact~\ref{Fact:GM} there is  one extension over $\L_p(\k)$ ramified at   $v$ \black  and thus it must be  $\L_p(\k)\tilde{\k}$.
Fact~\ref{Fact:Chevalley} applied to $\L_p(\k)\tilde{\k}/\L_p(\k)$ implies $(\#\Cl_{\L_p(\k)\tilde{\k}},p)=1$  so $\L_p(\tilde{\k})=\L_p(\k)\tilde{\k}$. Proposition~\ref{prop:stability} gives  $f(\tilde{\k}) > f(\k)$. 
\black
\end{proof}

\subsection{The inductive step}
\begin{Theorem}\label{Theorem:Induction}
%Let the hypotheses be as in Proposition~\ref{prop:Ozaki} except we 
Assume that $\L_p(\k)/\k$ has $\lambda_\k \geq 2h^1(\G)+3$ Minkowski units. Let $1\to \Zp\to \G'\to \G \to 1.$
\black If $\mu_p \not\subset \k$ (resp. $\mu_p \subset \k$) there exists a field $\k'/\k$ that is  a $\Zp$-extension (resp.  a compositum of two successive $\Zp$-extensions)  such that $\Gal(\L_p(\k')/\k') \simeq \G'$ and $\L_p(\k')/\k'$ has at least $2h^1(\G')+3$ Minkowski units.
\end{Theorem}
\begin{proof}
Recall that our finite tame primes split completely in $\K/\k$.

We  first treat the split case. This is independent of whether or not $\mu_p \subset \k$.

\underline{Split case.}
Choose tame degree one primes  $w_1$ and $w_2$ of $\Gov(\k)\K$  such that
\begin{itemize}
\item  $\Fr_{w_1}  =((x_1,x_2,\cdots, x_{h^1(\G)},0,\cdots,0),0) \in \Gal( \Gov(\L_p(\k))/\Gov(\k)\K) \subset M(\L_p(\k))$. This is possible as the tuple lies in $I_\G\cdot M(\L_p(\k))$ and $\Gov(\k)\K \subset F$. 
As $\Fr_{w_1}$ projects to $0$ in $\Gal(F/\K)$, we see
for  $v_1$ of $\k$ below $w_1$ that $\langle \Fr_{v_1}\rangle_{F/\K} =0$ so 
by Proposition~\ref{prop:Ozaki}
$\langle \Fr_{v_1}\rangle_{\Gov(\k)/\k(\mu_p)}=0$.
By Fact~\ref{Fact:GM} applied to $\k$ there is one $\Zp$-extension  $D_1/\k$ ramified at $v_1$. 
Fact~\ref{Fact:GM} also gives 
(see the proof of Theorem~\ref{Theorem:basechange} as well) a unique   $\Zp$-extension  of $\L_p(\k)$ ramified at $v_1$, namely $D_1\L_p(\k)/\L_p(\k)$.
  \item $\Fr_{w_2} =((0,0,\cdots ,0_{h^1(\G)},x_1,x_2,\cdots,x_{h^1(\G)}, 0,0,0,\cdots 0),0)   $ so for $v_2$ of $\k$ below $w_2$, $\langle \Fr_{v_2}\rangle_{F/\K} =0$.
 %that for some $w_2$ above $v_2$ in $\L_p(\k)$ we have  $$\langle \Fr_{w_2}\rangle  =  \langle ((0,0,\cdots ,0,x_1,x_2,\cdots,x_{h^1(\G)}, 0,0,0,\cdots 0),0)\rangle .$$  
 We also insist that  $v_2$ 
remains prime in $D_1/\k$. This last condition is linearly disjoint from the rest of the defining splitting conditions on $v_2$ and imposes no contradiction.  
 Again, there is one $\Zp$-extension  of both $\k$ and $\L_p(\k)$ ramified at $v_2$, namely $D_2/\k$.
Let $D/\k$ be a `diagonal' extension between $D_1$ and $D_2$   ramified at both $v_1$ and $v_2$.  There are $p-1$ of these.
\black
\end{itemize}
 
 Fact~\ref{Fact:GM} and our choices of the Frobenii of %the primes of $\L_p(\k)$ above 
 $v_1$ and $v_2$ imply  $h^1\big(\Gal(\L_p(\k)_{\{v_1,v_2\}}/\L_p(\k)\big)=2$. 
 %This  requires the existence of $2h^1(\G)$ Minkowski units in  $\L_p(\k)/\k$. 
%Had there been 
(With only  $h^1(\G)$ Minkowski units, we would have  had 
$h^1\big(\Gal(\L_p(\k)_{\{v_1\}}/\L_p(\k)\big)=h^1\big(\Gal(\L_p(\k)_{\{v_2\}}/\L_p(\k)\big)=1$, but  
$h^1\big(\Gal(\L_p(\k)_{\{v_1,v_2\}}/\L_p(\k)\big) >2$.)
$$ 
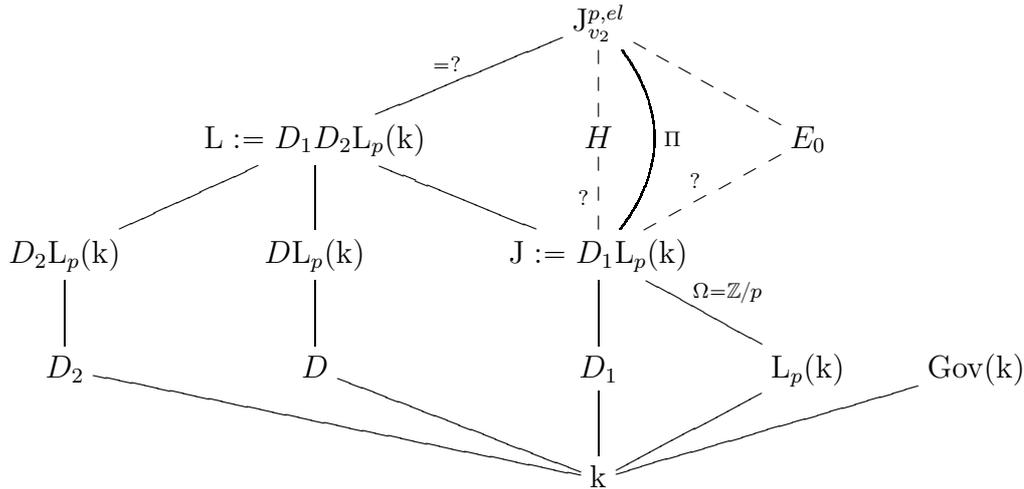
{
\xymatrix{ 
& &  \J_{v_2}^{p,el} \ar@{-}[dl]_{=?} \ar@{-}@/^1.75pc/[dd]^\Pi \\
&\L:=D_1D_2\L_p(\k)& H \ar@{--}[d]_? \ar@{--}[u]& E_0 \ar@{--}[dl]_? \ar@{--}[ul]&\\
D_2\L_p(\k)\ar@{-}[ur] & D\L_p(\k) \ar@{-}[u] & \J:=D_1\L_p(\k) \ar@{-}[ul] &  &\\
%& & & &\L_p(\k)\ar@{-}[ull]_{\Omega=\Zp}\\
D_2 \ar@{-}[u] & \ar@{-}[u] D & D_1  \ar@{-}[u] & \L_p(\k)\ar@{-}[ul]_{\Omega=\Zp}& \Gov(\k)\ar@{-}[lld]\\
&&  \k\ar@{-}\ar@{-}[u] \ar@{-}[ul] \ar@{-}[ull] \ar@{-}[ur] &  & &\\
}}$$

\smallskip

Set $\L:=D_1D_2\L_p(\k)$, $\J:=D_1\L_p(\k)$ and
note  $\L/D$ is unramified as $D/\k$ has absorbed all ramification at $\{v_1,v_2\}$. 
We will solve the problem by showing $(\# \Cl_{D_1D_2\L_p(\k)},p)=1$.

\smallskip
Since $(\# \Cl_{\L_p(\k)},p)=1$ and our choice of $v_1$ is such that  
$h^1(\Gal(\L_p(\k)_{\{v_1\}}/\L_p(\k))=1$,  Fact~\ref{Fact:Chevalley} applied to $\J/\L_p(k)$ implies $(\# \Cl_{\J},p)=1$.

\smallskip
We now prove that there exists a unique $\Z/p$-extension over $\J$ unramified outside $v_2$,  namely $\L$. 
Set $\Omega=\Gal(\J/\L_p(\k))$,  $\J_{\{v_2\}}^{p,el}$  to be the maximal elementary $p$-abelian extension of $\J$ inside 
$\J_{\{v_2\}}$, and $\Pi=\Gal(\J_{\{v_2\}}^{p,el}/\J)$.    Then $\Omega$ acts on $\Pi$  and trivially on $\Gal(\L /\J)$. We claim this is the only  $\Zp$-extension of $\J$ in $\J_{\{v_2\}}^{p,el}/\J$  on which $\Omega$ acts trivially: 
If not, there exists another $\Z/p$-extension $H/\J$ unramified outside $v_2$ and Galois over $\L_p(\k)$. Hence $\Gal(H/\L_p(\k))$ has order $p^2$ and is abelian. The extension $H/\L_p(\k)$ cannot be cyclic %of degree $p^2$ 
because  all inertia elements have order $p$ and would then fix an everywhere unramified extension of $\L_p(\k)$, a contradiction.
Suppose now that $\Gal(H/\L_p(\k))\simeq \Z/p \times \Z/p$, with $H \neq \J D_2=\L$. Then $\Gal(H D_2/\L_p(\k))\simeq (\Z/p)^3$: this contradicts the already established fact that $h^1(\Gal(\L_p(\k)_{ \{v_1,v_2\}}/\L_p(\k))=2$.

\smallskip The final possibility is that 
 there exists a $\Z/p$-extension $E_0/\J$ unramified outside $v_2$,  different from $\L/\J$ and not fixed by $\Omega$; let $S_0$ be the set of ramification of $E_0/\J$. As primes above $v_2$ in $\L_p(\k)$ are inert in $\J/\L_p(\k)$, $\Omega(S_0)=S_0$: 
 then $\Omega$ takes $E_0$ to another $\Z/p$-extension $E_1/\J$  exactly ramified at $S_0$ and such that   $E_1\neq E_0$. The compositum $E_1E_0/\J$ contains a $\Z/p$-extension $E_0'/\J $ exactly ramified at a set $S_0' \subsetneq S_0$. Observe  that $E_0'\neq \L$ since  $\L/\J$ is totally ramified at every prime above $v_2$. Continuing the process, we obtain an unramified  $\Z/p$-extension  $H/\J$, which is impossible since $(\#\Cl_\J,p)=1$.
Thus $\L/\J$ is the unique  $\Z/p$-extension  unramified outside $v_2$. Fact~\ref{Fact:Chevalley} applied to $\L/\J$  implies  $(\#\Cl_L,p)=1$.

\smallskip

We have solved the split embedding problem with $\k'=D$ and $\Gal(\L_p(\k')/\k') = \G \times \Zp$. It required one base change ramified at two tame finite primes. 
Proposition~\ref{prop:growth} implies $f(\k')\geq 2h^1(\G')+3$ so
the induction can proceed. 

\vskip1em

For the nonsplit case we treat  $\mu_p \not \subset \k$ and $\mu_p \subset \k$ separately. Theorem~\ref{Theorem:basechange} is only used in the nonsplit case when $\mu_p \subset \k$.

\vskip1em
\underline{The nonsplit case, $\mu_p \not \subset \k$.}
By Lemma~\ref{lemm:lemm_lifting} we may use one tame prime $v$ of $\k$ 
to find a {\it ramified} solution to the embedding problem. 
As $\mu_p \not \subset \k$ implies $\Gov(\k) \cap \L_p(\k) = \k$, we can assume $v$ splits completely in $\K/\k$.
%- this follows from Lemma~\ref{lemm:lemm_lifting}. 
%A prime $w$ of $\L_p(\k)$ above $v$ has 
Choosing any $w|v$ of $\K$ we set
$\Fr_{w} = ((z_1,z_2,\cdots, z_{\lambda_\k}),n_0) \in M(\L_p(\k))$
 %$$\langle  \Fr_w \rangle = \langle ((z_1,z_2,\cdots, z_\lambda),n) \rangle \subset \F_p[\G]^\lambda \oplus N$$ 
%in $\F_p[\G]^\lambda\oplus N$ \black 
where  we claim $n_0 \notin I_\G\cdot N$ 
 and $z_i \in I_\G \subset \F_p[\G]$.  Indeed, if any $z_i \notin I_\G$, its $\F_p[\G]$-span is all of $\F_p[\G]$ and
by Fact~\ref{Fact:GM}
there is no $\Zp$-extension of $\L_p(\k)$ ramified at the $w|v$, contradicting that we are solving an embedding problem with $v$.
If $n_0 \in I_\G \cdot N$,  then the projection of $\Fr_w$ to $\Gal(F/\K)$ is trivial so Proposition~\ref{prop:Ozaki} implies $\langle \Fr_v \rangle_{\Gov(\k)/\k(\mu_p)}=0$ and the embedding problem we are solving is split, also a contradiction. 

%As $\mu_p \not\subset \k$ one  sees that $\L_p(\k) \cap \Gov(\k)=\k$ so 
\smallskip
Choose a degree one  $w_1$ of  $\K$ with %splits completely in $\K/\k$ and 
%for some $w_1$ above $v_1$ in $\L_p(\k)$ we have 
 $ \Fr_{w_1} = ((x_1,x_2,\cdots,x_{h^1(\G)}, 0,0,0,\cdots 0),n_0) \in M(\L_p(\k)) $ where $n_0$ is as in the previous paragraph. \
  %where $\{x_1,x_2,\cdots,x_{h^1(\G)} \}$ generate $I_\G$, the augmentation ideal of $\G$.
  Let $v_1$ be the prime of $\k$ below $w_1$.
 By Fact~\ref{Fact:GM} (also see the proof of Theorem~\ref{Theorem:basechange}) there is one $\Zp$-extension  $D_1/\L_p(\k)$ ramified at $v_1$.
 
 \smallskip
Choose a  degree one \black $w_2$ of  $\K$ %(the prime of $\k$ beneath it is $v_2$) 
with %splits completely in $\K/\k$ and 
%that for some $w_2$ above $v_2$ in $\L_p(\k)$ we have 
 $\Fr_{w_2} =((0,0,\cdots ,0,x_1,x_2,\cdots,x_{h^1(\G)}, 0,0,0,\cdots 0),n_0)\in M(\L_p(\k))   $
 and the primes of $\L_p(\k)$ above $v_2$ remain prime in $D_1/\L_p(\k)$. 
 This last condition is linearly disjoint from the splitting conditions defining $v_2$ and imposes no contradiction. 
Again by Fact~\ref{Fact:GM} there is one $\Zp$-extension  $D_2/\L_p(\k)$ ramified at $v_2$.

\smallskip
As the free components of 
 of $\Fr_w$, $\Fr_{w_1}$ and $\Fr_{w_2}$  are all in $I^{\lambda_\k}_\G$, their projections 
 to $\Gal(F/\K)$ 
 depend only on $n_0$ and
% are all the same and
Proposition~\ref{prop:Ozaki}   implies
  $ 0 \neq \langle \Fr_{v} \rangle_{\Gov(\k)/\k} =\langle \Fr_{v_1}\rangle_{\Gov(\k)/\k}   =
  \langle \Fr_{v_2}\rangle_{\Gov(\k)/\k}$. 
Thus  there is no extension of $\k$
ramified at either $v_1$ or $v_2$, but, by Fact~\ref{Fact:GM}, there is a $\Zp$-extension of $\k$
 ramified at $\{v_1,v_2\}$. 
 \black Call it $D$. Note $\G'\simeq \Gal(D_1/\k)\simeq \Gal(D_2/\k) \simeq \Gal(D_1D_2/D).$
$$\black {
\xymatrix{ 
&D_1D_2=DD_1=DD_2 \ar@{-}[dr]&&\\
D_1\ar@{-}[ur]& D_2 \ar@{-}[u]& D\L_p(\k)&  &\\
%& \L_p(\k)\ar@{-}[ur] \ar@{-}[ul] \ar@{-}[u]& & &\\
\Gov(\k) \ar@{-}[dr] & \L_p(\k)\ar@{-}[ur] \ar@{-}[ul] \ar@{-}[u]& D\ar@{-}[u] \ar@{-}[u]_{\G}&  &\\
& \k  \ar@{-}[u]_{\G} \ar@{-}[ur]&  & &\\
}}$$

That $D_1D_2$ has trivial $p$-class group   follows exactly as it did in the split case \black and we may set $\k'=D$ 
so $\L_p(\k') = D_1D_2$ and $\Gal(\L_p(\k')/\k') \simeq \G'$.

\smallskip

We have solved the embedding problem in the nonsplit case when $\mu_p \not \subset \k$.
We performed one base change ramified at two tame finite primes and  
Proposition~\ref{prop:growth} implies 
$f(\k')\geq 2h^1(\G')+3$ so 
the induction can proceed.

\vskip1em\noindent
\underline{The nonsplit case, $\mu_p \subset \k$.}
%We finally address the case $\mu_p \subset \k$. 
We can no longer assume $\L_p(\k)\cap \Gov(\k) =\k$. 

Let $0\neq \varepsilon \in \sha^2_{\k,\emp}$ be the obstruction to our embedding problem $\G' \twoheadrightarrow \G$. Using Lemma~\ref{lemm:lemm_lifting}, 
let $v$ of $\k$ be a tame prime annihilating $\varepsilon$. The difficulty is that in the diagram 
below we may have   $\L_p(\k)\cap \Gov(\k) \supsetneq \k$ {\it and} that $\Fr_v$, which is necessarily nonzero
in $M(\k)$, may also be  nonzero in $\Gal((\L_p(\k)\cap \Gov(\k)) /\k)$. This prevents us from also choosing
$v$ to split completely  in $\L_p(\k)/\k$ and as we need in $ \Gov(\L_p(\k))/\L_p(\k)$  to ensure there is only one extension of $\L_p(\k)$ ramified at the primes of $\L_p(\k)$ above $v$.  If we could choose $v$ to annihilate $\varepsilon$ such that 
$\Fr_v =  0  \in \Gal(\L_p(\k)/\k)$, we would be able to proceed as in the  $\mu_p \not \subset k$ case. 
%to solve the embedding problem.  Henceforth we assume this is not possible. 
We  get around this by a base change.

\smallskip
By Kummer theory and the definition of governing fields,
 %and Kummer theory, for any number field~$\L$, 
 $\Gal(\Gov(\L)/\L(\mu_p))$ is an elementary $p$-abelian group. \black
Let $\tilde{\k}/\k$ be a tamely ramified $\Zp$-extension as given by Theorem~\ref{Theorem:basechange} so
$\Gal(\L_p(\tilde{\k})/\tilde{\k}) =\G$. By Proposition~\ref{prop:stability} we have $\lambda_{\tilde{\k}} \geq 2h^1(\G)+3$. 
%for $\L_p(\tilde{\k})/\tilde{\k}$.

$$\black 
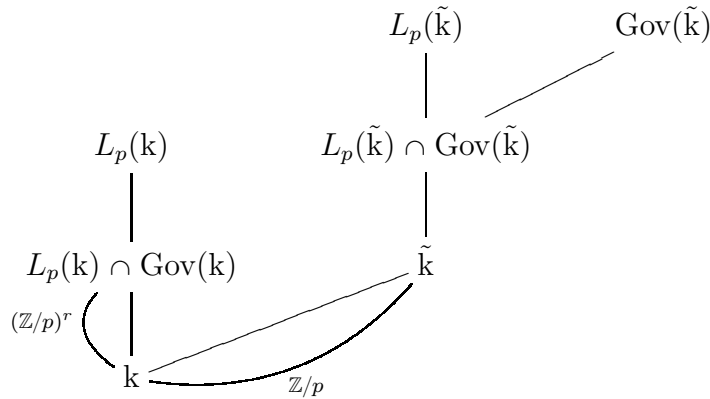
{
\xymatrix{ 
%& L_p( \tilde{\k}) \\
&  L_p( \tilde{\k}) & \Gov(\tilde{\k})\\
L_p(\k)& L_p(\tilde{\k}) \cap \Gov(\tilde{\k}) \ar@{-}[ur] \ar@{-}[u]& \\
L_p(\k) \cap \Gov(\k)\ar@{-}[u]& \tilde{\k}\ar@{-}[u]\\
\k\ar@{-}[u]\ar@{-}[ur]  \ar@{-}@/_1.5pc/[ur]_{\Zp} \ar@{-}@/^1.5pc/[u]^{(\Zp)^r} 
}}$$

As $\Gov(\k)\cap \tilde{\k} =\k$, we may choose a prime $v$ to solve the embedding problem for $\k$ whose Frobenius is nontrivial in $\Gal(\tilde{\k}/\k)$, that is $v$ remains prime in $\tilde{\k}/\k$. As observed above,
$\L_p(\tilde{\k}) \cap \Gov(\tilde{\k})/\tilde{\k}$ is a $(\Zp)^r$-extension for some $r$ and, as $\Gal(\L_p(\k)/\k) = 
\Gal(\L_p(\tilde{\k})/\tilde{\k}) =\G$, it is the base change of such  a subextension of $\L_p(\k)/\k$ from $\k$  so
$\L_p(\tilde{\k}) \cap \Gov(\tilde{\k})/\k$ is a $(\Zp)^{r+1}$-extension. Since $v$ remains prime in $\tilde{\k}/\k$ and residue field extensions are cyclic, it splits completely in $ L_p(\tilde{\k}) \cap \Gov(\tilde{\k})/\tilde{\k}$. As the embedding problem is solvable over $\k$ by allowing ramification at $v$, it is also solvable over $\tilde{\k}$ by allowing ramification at the unique prime of $\tilde{\k}$ above $v$. Thus 
$\varepsilon \in \sha^2_{\tilde{\k},\emp} \hookrightarrow \CyB_{\tilde{\k},\emp} =M(\tilde{\k})$
actually lies in $\Gal\left( \Gov(\tilde{\k})/ \left(\L_p(\tilde{\k}) \cap \Gov(\tilde{\k})\right)\right)$. The base change shifted  the obstruction to outside of our $p$-Hilbert class field tower! 
%We now proceed as in the $\mu_p \not \subset \k$ case. 
The rest of the proof is identical to the $\mu_p \not \subset \k$ case.
\end{proof}

We now prove the Main Theorem of the Introduction:
\begin{proof}
We have verified the base case of the induction in Proposition~\ref{prop:finish} and the inductive step with Theorem~\ref{Theorem:Induction}. It remains to count degrees and ramified primes. 
Proposition~\ref{prop:finish}  involved three $\Zp$-base changes, the first two ramified at one tame prime and the last at two tame primes.  
The inductive steps breaks into cases as follows
\begin{itemize}
\item $\mu_p \not \subset \k_0$: At each of the $\log_p(\#\HH) -1$ inductive stages we need one base change ramified at two primes for a total of $3+\left(\log_p(\#\HH)-1\right)$ base changes ramified at  $4+2(\log_p(\#\HH)-1)$ primes.
\item $\mu_p  \subset \k_0$: At each of the $\log_p(\#\HH) -1$ inductive stages we need two base changes, the first ramified at one prime and the second at two primes. There are  $3+2\left(\log_p(\#\HH)-1\right)$ base changes 
 ramified at  $4+3(\log_p(\#\HH)-1)$ primes.
\end{itemize}
\end{proof}

\black

\end{document}